\magnification\magstep1
\hsize 15.5truecm
\scrollmode
\overfullrule=0pt
\input amssym.def
\input amssym.tex
 at 17.28truept
 at 14.4truept
 at 12truept
 at 10.95truept
 at 17.28truept
 at 10truept

\def\title#1{%
\vskip0pt plus.3\vsize\penalty-100%
\vskip0pt plus-.3\vsize\bigskip\vskip\parskip%
\bigbreak\bigbreak\centerline{\bf #1}\bigskip%
}

\def\chapter#1#2{\vfill\eject
\centerline{\bf Chapter #1}
\vskip 6truept%
\centerline{\bf #2}%
\vskip 2 true cm}

\def\section#1#2{%
\def\\{#2}%
\vskip0pt plus.3\vsize\penalty-100%
\vskip0pt plus-.3\vsize\bigskip\vskip\parskip%
\par\noindent{\bf #1\hskip 6truept%
\ifx\empty\\{\relax}\else{\bf #2\smallskip}\fi}}

\def\subsection#1#2{%
\def\\{#2}%
\vskip0pt plus.3\vsize\penalty-20%
\vskip0pt plus-.3\vsize\medskip\vskip\parskip%
\def\TEST{#1}%
\noindent{\ifx\TEST\empty\relax\else\bf #1\hskip 6truept\fi%
\ifx\empty\\{\relax}\else{#2\smallskip}\fi}}

\def\proclaim#1{\medbreak\begingroup\noindent{\bf #1.---}\enspace\sl}

\def\endproclaim{\endgroup\par\medbreak}


\def\comfig#1#2\par{
\medskip
\centerline{\hbox{\hsize=10cm\eightpoint\baselineskip=10pt
\vbox{\noindent #1}}}\par\centerline{ Figure #2}}

\def\figcom#1#2\par{
\medskip
\centerline
{Figure #1}
\par\centerline{\hbox{\hsize=10cm\eightpoint\baselineskip=10pt
\vbox{\noindent #2}}}}
\def\bull{~\vrule height .9ex width .8ex depth -.1ex}


\def\comfig#1#2\par{
\medskip
\centerline{\hbox{\hsize=10cm\eightpoint\baselineskip=10pt
\vbox{\noindent{\sl  #1}}}}\par\centerline{{\bf Figure #2}}}

\def\figcom#1#2\par{
\medskip
\centerline
{{\bf Figure #1}}
\par\centerline{\hbox{\hsize=10cm\eightpoint\baselineskip=10pt
\vbox{\noindent{\sl  #2}}}}}

\def\em{\sl}

\def\\{\hfill\break}

\def\bull{~\vrule height .9ex width .8ex depth -.1ex}

\def\a{\alpha}
\def\b{\beta}
\def\g{\gamma}

\def\la{\lambda}
\def\La{\Lambda}

\def\u{\upsilon}

\def\CC{{\bf C}}

\def\ZZ{{\bf Z}}

\def\cO{{\cal O}}

\def\cW{{\cal W}}

\def\tr{\mathop{\rm tr}\limits}


\def\la{\lambda}

\def\eps{\varepsilon}

\def\sqr#1#2{{\vcenter{\hrule height.#2pt%
\hbox{\vrule width.#2pt height#1pt\kern#1pt%
\vrule width.#2pt}%
\hrule height.#2pt}}}

\newfam\gotfam
\font\twlgot=eufm10 at 12pt
\font\tengot=eufm10

\font\sevengot=eufm7
\textfont\gotfam=\twlgot
\scriptfont\gotfam=\tengot 
\scriptscriptfont\gotfam=\sevengot


\def\pr{\partial}
\def\res{{\rm res}}
\def\Ad{{\rm
Ad}}  \def\End{{\rm End}} \def\sh{{\rm sh}} 
  \def\w{\omega}\def\La{\Lambda}

\centerline{\bf Quantum groups in higher genus and Drinfeld's new
realizations method ($sl_{2}$ case).}
\bigskip
\centerline{B. Enriquez and V.N. Rubtsov}
\medskip
{\bf Abstract.} {\em We define double (central and cocentral) 
extensions of Manin
pairs introduced by Drinfeld, attached to curves and meromorphic
differentials. We define ``infinite twistings'' of these pairs and
quantize them in the $sl_{2}$ case, adapting Drinfeld's ``new
realizations'' technique. We study finite dimensional representations of
these algebras in level $0$, and some elliptic examples.}
\medskip

\section{}{Introduction.} 

In [5], V. Drinfeld introduced examples of Manin
pairs, attached to the data of a curve, a meromorphic differential on
it, and a finite dimensional reductive Lie algebra. He remarked, that
only in the cases where the curve had genus $\le 1$, could these Manin
pairs be given a structure of Manin triple; in these cases, the
quantization of these Manin triples gives rise to well-known Hopf
algebras (the Yangians, quantum
affine algebras and algebras connected with Sklyanin algebras). 
He raised the question of
quantizing these Manin pairs in the higher genus case, in the sense of
quasi-Hopf algebras. 

In this paper, we first present a double (central and cocentral)
extension of these Manin pairs. The general definition of these
extensions, in the case of Manin triples, is due to
M. Semenov-Tian-Shansky ([13]). This leads us to the problem of
the quantization of these extended Manin pairs. 

We then remark, that
this quantization problem can be approached in the spirit of the 
``new realizations'' of Drinfeld (introduced in [4] and developed
in [10], [3], [1]). This technique enabled Drinfeld to give a
quantum analogue of the passage from the Serre to the loop presentations
of an affine algebra; it can be presented as follows. The bialgebra
structure corresponding to quantum affine algebras, is a double
bialgebra structure. Conjugating the corresponding Manin triple by a
double group element, the bialgebra structure of the double gets
changed by a twisting (in the sense of [5]). 
Let us conjugate by affine Weyl group elements;
when their length tends to infinity, we get a limit Manin triple which
it is simple enough to quantize. The resulting Hopf algebra is then a
twisting of the one obtained by quantization of the initial Manin
triple (the Drinfeld-Jimbo Hopf algebra). 

In the present situation, we introduce a Lagrangian supplementary in
our Manin pair, and conjugate it by affine Weyl group elements as
before. (We note, that a family of supplementaries is provided by
a covering of the space of principal $G$-bundles over the curve $X$; 
we hint at a possible connection between the closedness of a $1$-form on it,
defined in terms of twisting, and  
a generalized classical Yang-Baxter identity, underlying
the integrability of the Hitchin system. 
We hope to return to this
question in [6].) 
In the limit, we obtain a Manin triple, whose quantization (in
the $sl_{2}$ case) is the main goal of this paper. Let us describe more
precisely its contents. 

Let $X$ be a smooth compact complex curve, $\w$ a meromorphic nonzero
one-form on $X$, $\{x_{i}\}\subset X$ the set of its zeroes and poles. 
Let for each $i$, $k_{x_{i}}$ be the local field
at $x_{i}$, $\cO_{x_{i}}$ the local ring at this point, $R\subset
\oplus_{i}k_{x_{i}}$
the ring of functions regular outside $\{x_{i}\}$. We choose a
supplementary $\La$ to $R$ in $\oplus_{i}k_{x_{i}}$, Lagrangian for the
scalar product defined by $\w$. To define a quantization of our Manin
triple, we need operators $A:R\to \oplus_{i}k_{x_{i}}$ and $B:\La\to 
\oplus_{i}k_{x_{i}}$, which serve to define $h-e$ and $h-f$ relations
(by $h-e$ relations, we understand relations between Fourier modes of
the quantum analogues of fields $h(z)$ and $e(z)$, etc.;
with $e,h,f$ the Chevalley generators of $sl_{2}$). 
These operators also provide us with $e-e$ and $f-f$ relations, which
appear in the form
$$
e(z)e(w)=a(z,w)e(w)e(z), \quad w\ll z; 
$$
Our aim is to put these relations in the form 
$$
(z-w+\sum_{i\ge 1}\hbar^{i}\a_{i}(z,w))e(z)e(w)
=(z-w+\sum_{i\ge 1}\hbar^{i}\b_{i}(z,w))e(w)e(z),
$$
$\a_{i},\b_{i}$ formal series in $z$, $w$, ($\hbar$ is the quantization
parameter) similar to the quantum affine
algebra relations
$$
(qz-w)e(z)e(w)=(z-qw)e(w)e(z).
$$
It turns out that to achieve this task, essentially one possibility for the 
operators $A$ and $B$ remains. The proof that it indeed leads to
$e-e$ and $f-f$ relations of the desired form, relies on a statement
about derivatives of a Green function (prop. 1), which allows to give
a universal treatment for all pairs $(X,\omega)$. The formal series
$\a_{i}$ and $\b_{i}$ are then obtained from formal solutions to certain
differential equations (eqs. (3.7)), where the variable is $\hbar$. 

The quantization we propose depends both on a choice of $\La$, and of a
certain element $\tau\in R\otimes R$. We show that the various
quantizations obtained are related to each other by twisting operations
(in the sense of [5]). 

We turn then to the problem of finite dimensional representations of our
algebras, at level $0$; these representations are indexed by points of formal
discs. We construct a family of $2$-dimensional representations. We
expect that higher spin representations can be constructed as well and
their tensor products have properties similar to those explained in
[2]. 

We close the paper by making explicit some examples. In a certain
elliptic case, we recover as $e-e$ relations some elliptic
$\cW$-algebra relations discovered in [8]. We also propose a twisting
of the ``automorphic'' Manin triples of [12], which could lead to a
``quantum currents'' description of the Hopf algebra arising in [15]. 

Further problems related to the present construction could be the
following: applying the Reshetikhin-Semenov method for constructing
central elements at the critical level ([11]); construction of level $1$
representations as in [9], vertex operators 
and the corresponding quantum Knizhnik-Zamolodchikov (KZ)
equations; generalization from $sl_{2}$ to an arbitrary semisimple Lie
algebra. The resulting quantum KZ equations at the critical level, might
then be considered as $q$-deformations of the holonomic systems of
equations occurring in the geometric Langlands program; working out such
$q$-deformations using quantum affine algebras was proposed by E. Frenkel
and N. Reshetikhin. Finally, 
in [7], P. Etingof and D. Kazhdan showed how to attach to any
associator, a quantization procedure for bialgebras. 
It would be interesting to understand, whether the construction
presented here can be obtained from the KZ
associator, as it is the case for finite
dimensional Lie algebras. 

We would like to thank B. Feigin, G. Felder and
E. Vasserot for discussions related to the subject of this paper, and
S. Khoroshkin for explaining to us the ideas of [4]. 
The work of V.R. was supported by the CNRS and partially by grant RFFI
95-01-01101. He expresses his thanks to the Centre de Math\'ematiques
de l'Ecole Polytechnique, for excellent working conditions.  

\section{}{1. Manin triple}

\it 1. Drinfeld's Manin pairs. \rm 

Let $X$ be a smooth compact complex curve, $\w$ a meromorphic nonzero
one-form on $X$, $\{x_{i}\}\subset X$ the set of its zeroes and poles. 
Let for each $i$, $k_{x_{i}}$ be the local field
at $x_{i}$, $\cO_{x_{i}}$ the local ring at this point, $R\subset
\oplus_{i}k_{x_{i}}$
the ring of functions regular outside $\{x_{i}\}$. Let ${\bf a}$ be a simple
complex Lie
algebra, $\langle,\rangle_{{\bf a}}$ its Killing form. Recall the Manin pair
defined by Drinfeld in [5]: 
endow 
$$
{\bf g}_{0}={\bf a}(\oplus_{i}k_{x_{i}})
$$ 
with
the bilinear
form $\langle x_{1}, x_{2}\rangle_{0}=\sum_{i}\res_{x_{i}}(\langle
x_{1},x_{2}\rangle_{{\bf a}}\omega)$; then 
${\bf a}(R)$ is a Lagrangian subalgebra of ${\bf g}_{0}$; this defines 
Drinfeld's
Manin pair
$$
({\bf a}(\oplus_{i}k_{x_{i}}), {\bf a}(R)). 
$$

\it 2. Double extension. \rm 

We extend this Manin pair
in the following way. Let $\partial$ be the derivation of 
${\bf a}(\oplus_{i}k_{x_{i}})$
defined by $\partial f=df/\omega$. We denote in the same way the
derivation of $\oplus_{i}k_{x_{i}}$, defined by the same formula. 
Let ${\bf g}$ be the
skew product of ${\bf g}_{0}$ by $\partial$; we have 
$$
{\bf g}={\bf g}_{0}\oplus\CC \tilde D,
$$
${\bf g}_{0}\subset {\bf g}$ is a Lie algebra 
homomorphism, and $[\tilde D,x]=\partial x$
for $x\in {\bf g}_{0}$. Since $\partial$ preserves ${\bf a}(R)$,
${\bf a}(R)\oplus\CC\tilde D$ is a Lie subalgebra of ${\bf g}$. 
Let $\hat{\bf g}$ be the central extension of ${\bf g}$ by $\CC
K$ using the cocycle defined by $c(x,y)=\res_{x}\langle x,dy\rangle_{a}K$,
$c(x,\tilde D)=0$ for
$x\in {\bf g}_{0}$. We have 
$$
\hat{\bf g}={\bf g}\oplus \CC K,
$$
with the usual commutation rules. 
Since $c$ vanishes on ${\bf a}(R)\oplus\CC\tilde D$, this
algebra has a section to $\hat{\bf g}$, that we denote by
${\bf g}_{R}$. 
Identifying $\hat{\bf g}$ with ${\bf g}\oplus\CC K$, 
${\bf g}_{R}$ is identified  with
$$
({\bf a}(R)\oplus\CC\tilde D)\times\{ 0\}.
$$
Let $D=(\tilde D,0)$. 
Let us consider now on $\hat{\bf g}$ the symmetric bilinear form,
defined by 
$\langle K, D\rangle=2$, $\langle D, {\bf g}\times \{0\}\rangle=0$, $\langle
K, {\bf g}_{0}\times \{0\}\rangle=0$, $\langle (x_{1},0),(x_{2},0)\rangle=
\langle x_{1},x_{2}\rangle_{0}$ for $x_{1},x_{2}\in {\bf g}_{0}$. Then
$\langle ,\rangle$ is invariant, and ${\bf g}_{R}$ is a subalgebra of
$\hat{\bf g}$, Lagrangian w.r.t. this scalar product. 
$$
(\hat{\bf g},{\bf g}_{R})
$$ 
is a
double extension (central and cocentral) of the above Manin pair. 

\it 3. Lagrangian supplementaries. \rm 

Consider on $\oplus_{i}k_{x_{i}}$, the
scalar product defined by $\langle f_{1},
f_{2}\rangle_{\oplus_{i}k_{x_{i}}}=\sum_{i}
\res_{x_{i}}(f_{1}f_{2}\w)$. $R$ is a subspace of $\oplus_{i}k_{x_{i}}$,
Lagrangian
w.r.t. this scalar product. Fix a Lagrangian supplementary $\La$ to $R$,
commensurable with $\oplus_{i}\cO_{x_{i}}$.
Then 
$$
{\bf a}\otimes \La\oplus\CC K
$$ 
is a Lagrangian supplementary to ${\bf g}_{R}$ in $\hat{\bf g}$. 

We note here, that a family of such supplementaries can be defined
in the following way. Let $\La_{0}$ be a Lagrangian subspace of
$\oplus_{i}k_{x_{i}}$, containing $\oplus_{i}\cO_{x_{i}}$. Let $G$ be a
Lie group, with Lie algebra ${\bf a}$. For $g\in
G(\oplus_{i}k_{x_{i}})$, let 
$$
^{g}{\bf a}(\La_{0})=
g {\bf a}(\La_{0})g^{-1}; 
$$
for generic $g$ this defines a Lagrangian supplementary to ${\bf g}_{R}$
in $\hat{\bf g}$, which up to equivalence depends only on the class of
$g$ in $G(R)\setminus G(\oplus_{i}k_{x_{i}})/{\rm Stab}\ {\bf
a}(\La_{0})$. All the resulting bialgebra structures on ${\bf g}_{R}$
are then associated by twisting. 
The twisting between two bialgebra structures associated to nearby
points defines an element of
$\wedge^{2}{\bf a}(R)$; we thus get a $1$-form $\Phi\in
\Omega^{1}(G(\oplus_{i}k_{x_{i}}),
\wedge^{2}{\bf a}(R))$, equivariant
w.r.t. left $G(R)$-translations. This $1$-form is closed; we hope that 
the expression
of this fact can be interpreted as the generalized Yang-Baxter
identity for the dynamical $r$-matrices of the Hitchin system ([6]).

\it 4. Infinite twisting. \rm 

Let us return to the triple formed by $\hat{\bf g}$, ${\bf g}_{R}$ and
the Lagrangian 
supplementary ${\bf a}\otimes \La\oplus\CC K$, and conjugate it by affine Weyl
group elements. In the limit, we obtain the Manin triple 
$$
(\hat{\bf g}, {\bf g}_{+}, {\bf g}_{-})
$$
with 
$$
{\bf g}_{+}=
{\bf h}(R)\oplus {\bf n}_{+}(\oplus_{i}k_{x_{i}})\oplus\CC D, 
\quad {\bf g}_{-}={\bf h}\otimes\La\oplus
{\bf n}_{-}(\oplus_{i}k_{x_{i}})\oplus\CC K.
$$
Here ${\bf h}$, ${\bf n_{+}}$ and ${\bf n_{-}}$ are the
Cartan and opposite
nilpotent subalgebras of ${\bf a}$. Our aim will be to give a quantization of
the Lie bialgebras $({\bf g}_{\pm},\delta_{\pm})$, where $\delta_{\pm}$ are
the cobrackets 
defined by the above Manin triple, in the case where ${\bf a}=sl_{2}(\CC)$. 

${\bf g}_{+}$ and ${\bf g}_{-}$ can be presented as follows. 
${\bf g}_{+}$ has generators
$D$, $h^{+}(r)$, $r\in R$, and
$e(\eps)$, $\eps\in \oplus_{i}k_{x_{i}}$, with
$$
h^{+}(\a_{1}r_{1}+\a_{2}r_{2})=\a_{1}h^{+}(r_{1})+\a_{2}h^{+}(r_{2}),
\quad\a_{i}\in\CC, r_{i}\in R,
\leqno(1.1)
$$
and
$$
e(\a_{1}\eps_{1}+\a_{2}\eps_{2})=\a_{1}e(\eps_{1})+\a_{2}e(\eps_{2}),
\quad
\a_{i}\in\CC, \eps_{i}\in \oplus_{i}k_{x_{i}};\leqno(1.2)
$$
${\bf g}_{-}$ has generators
$K$, $h^{-}(\la)$, $\la\in \La$, and
$f(\eps)$, $\eps\in \oplus_{i}k_{x_{i}}$, with
$$
h^{-}(\a_{1}\la_{1}+\a_{2}\la_{2})=\a_{1}h^{-}(\la_{1})+\a_{2}h^{-}(\la_{2}),
\quad
\a_{i}\in\CC, \la_{i}\in \La,\leqno(1.3)
$$
$$
f(\a_{1}\eps_{1}+\a_{2}\eps_{2})
=\a_{1}f(\eps_{1})+\a_{2}f(\eps_{2}),\quad
\a_{i}\in\CC, \eps_{i}\in \oplus_{i}k_{x_{i}}.\leqno(1.4)
$$
The relations are the
following. 
Let $e^{i}$, $e_{i}$ be dual bases of $R$ and $\La$ and $\sum
\eps^{i}\otimes \eps_{i}=\sum e^{i}\otimes e_{i}+e_{i}\otimes e^{i}$. 
Let $z=(z_{i})$ be a system of local coordinates at each point $x_{i}$. 
We 
define the formal series
$$
\eqalign{e(z)=\sum e(\eps^{i})\eps_{i}(z), \quad \quad f(z)=\sum
f(\eps^{i})\eps_{i}(z), 
\cr
h^{+}(z)=\sum h^{+}(e^{i})e_{i}(z), \quad
h^{-}(z)=\sum h^{-}(e_{i})e^{i}(z);
}\leqno(1.5) 
$$
then the Lie algebra relations for ${\bf g}_{+}$ and ${\bf g}_{-}$ are 
respectively  
$$
\eqalign{
[h^{+}(r), h^{+}(r')] =0,\quad [h^{+}(r),e(z)]=2r(z)e(z),
\quad r,r'\in R\cr
[D,h^{+}(r)] =h^{+}(\partial r), \quad
[D,e(z)]=-\pr_{z}e(z),
\quad
[e(z),e(w)]=0 \cr
}\leqno(1.6)
$$
and
$$
\eqalign{
[h^{-}(\la), h^{-}(\la')]=2\res_{x}(\la d\la')K,
\quad [h^{-}(\la),f(z)]=-2{\la(z)}
f(z),\cr [f(z),f(w)]=0,\quad
[K,h^{-}(\la)]=[K,f(\la')]=0, \quad \la,\la'\in\La\cr}
\leqno(1.7)
$$
Equalities (1.6) and (1.7) give 
$$
[h^{+}(z),e(w)]=2a_{0}(w,z)e(w), \quad 
[h^{-}(z),e(w)]=2a_{0}(z,w)e(w);
$$
we set $a_{0}(z,w)=\sum_{i}e^{i}(z)e_{i}(w)$, $e^{i}$, $e_{i}$ being
dual bases of $R$ and $\La$. We have 
$$
a_{0}(z,w)=[{1\over r_{0}(w)}
\sum_{i\in\ZZ}(z/w)^{i}]_{z\to R}\leqno(1.8)
$$
(the index $z\to R$ means acting by $\pi\otimes 1$, $\pi$ the projection
of $\oplus_{i}k_{x_{i}}$ to $R$ parallel to $\La$, in a completion of
$k\hat
\otimes
k=\CC((z,w))$)
and 
$$
a_{0}(z,w)+a_{0}(w,z)= {1\over r_{0}(z)}\sum_{i\in\ZZ}(z/w)^{i}, \leqno(1.9)
$$
which can be proved viewing the l.h.s. as a kernel. 
Thanks to (1.9), we
rewrite equalities (1.6.b,c), and (1.7.a,b) as 
$$
\eqalign{[D, h^{+}(z)]=
(-\pr_{z}h^{+}(z))_{z\to\La}
\quad  
[h^{-}(z), h^{-}(w)]=\bar \gamma(z,w) K, \cr
[h^{+}(z), e(w)]=2a_{0}(w,z)e(w), \quad 
[h^{-}(z), f(w)]=-2a_{0}(z,w)f(w),\cr} \leqno(1.10)
$$ 
where 
$$
\bar
\gamma(z,w)=-(\partial_{z}+\partial_{w})a_{0}(z,w)\in\wedge^{2}R.\leqno(6)
$$ 
The Jacobi identity is ensured by
$((\partial_{z}+\partial_{w})a_{0}(z,w))_{z\to\La}=0$. This identity is a
consequence of (1.11), which is
proven as follows. 
Pose $\partial
e^{i}=\sum_{j}a^{i}_{j}e^{j}$, $\partial
e_{i}=\sum_{j}c^{j}_{i}e_{j}+d_{ij}e^{j}$; $a^{i}_{j}+c^{i}_{j}=0$ and
$d_{ij}+d_{ji}=0$ because $\partial$ is anti-self-adjoint. 
$-(\partial_{z}+\partial_{w})a_{0}(z,w)$ is then equal to
$\sum_{i,j}d_{ij}e^{i}\otimes e^{j}$, and so belongs to $R\otimes R$ and
is antisymmetric; this proves (1.11). 
The pairing between ${\bf g}_{+}$ and ${\bf g}_{-}$ is given by 
$$
\langle D, K\rangle=2, \quad \langle e(z), f(w)\rangle={1\over {r_{0}(z)}}
\sum_{i\in\ZZ}(z/w)^{i}, \quad \langle h^{+}(z),
h^{-}(w)\rangle=2a_{0}(w,z) \leqno(1.12) 
$$
The formulas for the cobracket of ${\bf g}_{+}$ and ${\bf g}_{-}$ are then
respectively 
$$\eqalign{
\delta_{+}(e(z))
=e(z)\wedge h^{+}(z), 
\quad \delta_{+}(h_{n}^{+})=0 \cr
\delta_{+}(D)=\sum_{i,j}\res_{z=x_{i}}\res_{w=x_{j}}
\gamma(z,w)(h^{+}(z)\wedge
h^{+}(w))\omega_{z}\omega_{w}
\cr}
\leqno(1.13)
$$
and
$$\eqalign{
\delta_{-}
(f(z))=h^{-}(z)\wedge f(z)
+K\wedge \pr_{z}f(z),
\cr 
\delta_{-}(h^{-}(\xi))= K\wedge
\pr_{z}h^{-}(z), 
\quad \delta_{-}(K)=0
\cr}\leqno(1.14)
$$
The fact (1.13) defines a cocycle can be checked directly using the
identities (1.9) and
$$
\sum_{i,j}\res_{z'=x_{i}}\res_{w'=x_{j}}(a(z',z)\gamma(z',w')h^{+}(w')
\omega_{z'}\omega_{w'})
=(\pr_{z}h^{+}(z))_{z\to R},\leqno(1.15)
$$ 
(as before, the index $z\to R$ has the meaning of applying $\pi$), 
which is equivalent to 
$\sum_{i}\res_{w'=x_{i}}(\gamma(z,\w')h^{+}(w')\omega_{w'})
=(\pr_{z}h^{+}(z))_{z\to R}$; since
$h^{+}(z)$ belongs to a completion of $\hat{\bf g}\otimes\La$, 
it is enough to check
this identity replacing $h^{+}$ by any element of $\La$, and in
this case if follows from $(\partial e_{i})_{R}=\sum_{j}d_{ij}e^{j}$. 

Let us describe how a change of $\La$ affects $a_{0}(z,w)\in
R_{z}((w))$. Let $\La'$ be another Lagrangian supplementary to $R$,
commensurable with $\cO_{x}$; let $\pi'$ be the projection of
$\oplus_{i}k_{x_{i}}$
onto $R$ parallel to $\La'$, and let $a_{0;\La'}(z,w)=(1\otimes \pi')
({1\over
{r_{0}(w)}}[\sum_{i\in\ZZ}(z/w)^{i}])$, $\gamma_{\La'}(z,w)=
-(\partial_{z}+\partial_{w})a_{0;\La'}(z,w)$. The projection
of $\La$ on $R$ parallel to $\La'$ gives some antisymmetric element
$r_{1}\in R\otimes R$, and to any such element corresponds such a
$\La'$. We have then 
$$
a_{0;\La'}(z,w)=a(z,w)-r_{1}(z,w), 
\quad
\gamma_{\La'}(z,w)=\gamma(z,w)+(\partial_{z}+\partial_{w})
r_{1}(z,w). \leqno(1.17)
$$

\section{}{2. Quantization of ${\bf g}_{+}$ and ${\bf g}_{-}$}

Let $\hbar$ be a formal parameter and $T$ be the operator
${{\sh\hbar\pr}\over{\hbar\pr}}:\oplus_{i }k_{x_{i}}\to
(\oplus_{i}k_{x_{i}})[[\hbar]]$; 
since $T$ is symmetric for $\langle , \rangle_{\oplus_{i}k_{x_{i}}}$,
the expression
$$
\sum_{i}Te^{i}\otimes e_{i}-e^{i}\otimes Te_{i}
$$
belongs to $S^{2}R$. Let us fix $\tau\in R\otimes R[[\hbar]]$, 
such that 
$$
\tau+\tilde\tau=\sum_{i}Te^{i}\otimes e_{i}-e^{i}\otimes Te_{i},
\leqno(2.1)
$$
where we denote $\tilde f(z,w)=f(w,z)$. Let $U$ be the operator 
from $\La$ to $R[[\hbar]]$, such that 
$$
\tau=\sum_{i}Ue_{i}\otimes e^{i};
$$
$U$ verifies 
$$
\sum (T+U)e_{i}\otimes e^{i}+e^{i}\otimes (T+U)e_{i}=\sum
Te^{i}\otimes e_{i}+Te_{i}\otimes e^{i}.
$$
Let $U_{\hbar}{\bf g}_{+}^{(0)}$  (resp. $U_{\hbar}{\bf g}_{-}$) 
be the algebra with generators $h^{+}(r)$,
$r\in R$, $e(\eps)$, $\eps\in\oplus_{i}k_{x_{i}}$, 
(resp. $K$, $h^{-}(\la)$,
$\la\in \La$, $f(\eps)$, $\eps\in\oplus_{i}k_{x_{i}}$), subject
to relations (1.1-4), organized in generating
series (1.5), and subject to the relations
$$
[h^{+}(r), e(w)]=2r(w)e(w), e(z)e(w)=e^{2\hbar\sum
((T+U)e_{i})(z)e^{i}(w)}e(w)e(z),\leqno(2.2)
$$
$$
\Delta_{+}
(h^{+}(r))=h^{+}(r)\otimes 1+1\otimes h^{+}(r), \Delta_{+}
(e(z))=e(z)\otimes \exp({\hbar}((T+U)h^{+}(z))
+1\otimes e(z);\leqno(2.3)
$$
and
$$
K {\rm \ is \ central,\  }
[h^{-}(z), f(w)]=-2
q^{K(\pr_{z}+\pr_{w})}\{\sum e^{i}(z)(T+U)e_{i}(w)
\}f(w),\leqno(2.4)
$$
$$
{\rm or\ }[h^{-}(\la), f(w)]=-2q^{K\pr}((T+U)(q^{-K\pr}\la)_{\La})(w)f(w), 
\leqno(2.5)
$$
$$
[h^{-}(z),h^{-}(w)]={2\over\hbar}(q^{K(\pr_{z}+\pr_{w})}
-q^{-K(\pr_{z}+\pr_{w})})\sum e^{i}(z)(T+U)e_{i}(w)\leqno(2.6)
$$
$$
f(z)f(w)=q^{K(\pr_{z}+\pr_{w})}\{
e^{2\hbar\sum e^{i}(z)((T+U)e_{i})(w)}
\}f(w)f(z),\leqno(2.7)
$$
$$\Delta_{-}(K)=K\otimes 1 + 1 \otimes K,
\Delta_{-}(h^{-}(\la))=h^{-}((q^{K_{2}\pr}\la)_{\La})\otimes 1+
1\otimes h^{-}((q^{-K_{1}\pr}\la)_{\La}),\leqno(2.8)
$$
$$
\Delta_{-}(f(z))=(q^{-K_{2}\pr}f)(z)\otimes 
\exp(-{\hbar} (q^{K_{1}\pr}h^{-})(z))
+1\otimes (q^{K_{1}\pr}f)(z),\leqno(2.9)
$$
where, thanks to the formula for $\Delta_{-}(K)$, we view $\Delta_{-}$
as a system of maps from $(U_{\hbar}{\bf g}_{-})_{K_{1}+K_{2}}$ to 
$(U_{\hbar}{\bf g}_{-})_{K_{1}}\otimes (U_{\hbar}{\bf g}_{-})_{K_{2}}$,
for variable scalars $K_{i}$, where
$(U_{\hbar}{\bf g}_{-})_{k}=U_{\hbar}{\bf g}_{-}/(K-k)$ for any
scalar $k$. 
(2.8.b) can also be written 
$$
\Delta_{-}
h^{-}(z)=(q^{-K_{2}\pr}h^{-})(z)\otimes 1 +1 \otimes
(q^{K_{1}\pr}h^{-})(z).
\leqno(2.10)
$$
We can also write the $h^{-}-h^{-}$ commutator as 
$$
\eqalign{
[h^{-}(z),h^{-}(w)] & ={1\over \hbar}(T_{z}+T_{w})(q^{K(\pr_{z}+\pr_{w})}-
q^{-K(\pr_{z}+\pr_{w})})a_{0}(z,w) \cr
& +{2\over\hbar}(q^{K(\pr_{z}+\pr_{w})}-
q^{-K(\pr_{z}+\pr_{w})})\sum_{i}e^{i}\otimes Ue_{i}; 
\cr}
\leqno(2.11)
$$
recall that $a_{0}(z,w)=\sum_{i}e^{i}(z)e_{i}(w)$ and that
$(\pr_{z}+\pr_{w})^{\ge 1}a_{0}(z,w)\in R_{z}\otimes R_{w}$; this shows
that the r.h.s. of the formula for $[h^{-}(z),h^{-}(w)]$ belongs to
$\wedge^{2}R$, as it should be.  
(We will show later on how to put the $e-e$ and $f-f$ relations in a
correct form.)

The skew antipodes have the form 
$$
S'_{+}(h^{+}(r))=-h^{+}(r),
S'_{+}(e(z))=-\exp(-\hbar ((T+U)h^{+}(z))e(z),
\leqno(2.12)
$$
$$
S'_{-}(K)=-K, 
S'_{-}(h^{-}(\la))=-h^{-}(\la), S'_{-}(f(z))=-\exp(\hbar
h^{-}(z))f(z).
\leqno(2.13)
$$
The pairing 
$$
\langle e(z), f(w)\rangle=\delta(z/w), \quad \langle h^{+}(r),
h^{-}(\la)\rangle={2\over \hbar}\langle r,\la\rangle_{\oplus_{i}k_{x_{i}}}
\leqno(2.14)
$$ 
(with $\delta(z/w)=\sum_{i\in\ZZ}(z/w)^{i}$)
extends to
a Hopf algebra pairing between $U_{\hbar}{\bf g}_{+}^{(0)}$ and 
$U_{\hbar}{\bf g}_{-}$. 

The double of $U_{\hbar}{\bf g}_{-}$ has then the additional relations 
$$
[h^{+}(r), h^{-}(\la)]={2\over\hbar}\langle (q^{K\pr}-q^{-K\pr})r,\la
\rangle_{k_{x}},
\leqno(2.15) 
$$
$$
[h^{+}(r),f(z)]=-2(q^{K\pr}r)(z)f(z),\leqno(2.16)  
$$
$$
[e(z),f(w)]=(q^{K\pr_{w}}\delta(z/w))q^{(T+U)h^{+}(z)}
-(q^{-K\pr_{w}}\delta(z/w))q^{-h^{-}(w)},\leqno(2.17) 
$$
$$
[h^{-}(\la),e(w)]=2[(T+U)((q^{K\pr}\la)_{\La})](w)e(w).\leqno(2.18) 
$$
Let us define 
$$
\tilde D(h^{+}(r))=h^{+}(\pr r),\leqno(2.19)
$$
$$
 \tilde
D(e(z))=-\pr_{z}e(z)+{\hbar\over 2}
[\pr(T+U)h^{+}-(T+U)(\pr h^{+})_{\La}](z)e(z)
\leqno(2.20)$$
where
$(\pr h^{+})_{\La}(z)=\sum h^{+}(e^{i})(\pr e_{i})_{\La}$, and the index
$\La$ denotes the projection of $\oplus_{i}k_{x_{i}}$ on $\La$ parallel
to $R$; 
$\tilde D$ extends to a derivation of $U_{\hbar}{\bf g}_{+}^{(0)}$. Let 
$U_{\hbar}{\bf g}_{+}$ be the algebra generated by $U_{\hbar}{\bf
g}_{+}^{(0)}$ and the element $D$, such that 
$$
[D,x]=\tilde D(x)\leqno(2.21)
$$
for $x\in U_{\hbar}
{\bf g}_{+}^{(0)}$. Let us extend $\Delta_{+}$ to $U_{\hbar}{\bf g}_{+}$ by
$$
\Delta_{+}(D)=D\otimes 1+1\otimes D-{\hbar\over 4}
\big\{ h^{+}[((T+U)e_{i})_{R}]\otimes h^{+}(\pr e^{i})+
h^{+}[(\pr(T+U)e_{i})_{R}]\otimes h^{+}(e^{i})\big\}.
\leqno(2.22) 
$$
Here the index $R$ denotes the projection of $\oplus_{i}k_{x_{i}}$ on
$R$ parallel to $\La$. 
Then $\Delta_{+}$ defines a Hopf algebra structure on
$U_{\hbar}{\bf g}_{+}$, dual to $U_{\hbar}{\bf g}_{-}$ if we extend the
pairing (2.14) by 
$$
\langle D,K\rangle=1/2, \quad \langle D,
h^{-}(\la)\rangle=\langle D, f(z)\rangle=0;
$$
the counit and skew antipode extend then to $\epsilon(D)=0$, 
$$
S'_{+}(D)=-D-{\hbar\over 4}\sum_{i}
\big\{
 h^{+}(\pr e^{i})h^{+}[((T+U)e_{i})_{R}]
+h^{+}(e^{i})h^{+}[(\pr(T+U)e_{i})_{R}]
\big\},
\leqno(2.23) 
$$
and the quantum double relations are 
$$
\eqalign{
[ D, h^{-}(z) ] =-\pr h^{-}(z) & + q^{-K\pr} \big\{
\pr[(T+U^{*})h^{+}]_{R}
-[(T\pr h^{+})_{R}+U^{*}(\pr h^{+})_{\La}]\big\}(z)
\cr
&
+{1\over
2}(q^{K\pr}-q^{-K\pr})[\pr U^{*} h^{+}-(\pr(Th^{+})_{\La})_{R}-U^{*}(\pr
h^{+})_{\La}](z),
\cr
}\leqno(2.24)
$$
where $U^{*}:\La\to R[[\hbar]]$ is the map dual to $U:\La\to
R[[\hbar]]$, and 
$$
[ D , f(z) ] =-\pr f(z)+{\hbar\over 2}q^{K\pr}[\pr(T+U)h^{+}-(T+U)(\pr
h^{+})_{\La}](z)f(z). 
\leqno(2.25)
$$

\section{3. $e-e$ and $f-f$ relations.}{}

\medskip
\noindent {\it 3.1. Construction of $\gamma$.}

Let us compute the endomorphism of $R$, defined by
$$
\rho(f)(z)=\res_{w}a_{0}^{2}(z,w)f(w)\omega_{w}.\leqno(3.1) 
$$
$\rho$ is the restriction to $R$ of the endomorphism $\bar\rho$ of
$k_{x}$, defined by
$$
\bar\rho(f)=\res_{w}(a_{0}^{2}-\tilde
a_{0}^{2})(z,w)f(w)\omega_{w}.\leqno(3.2) 
$$
Let
$\a(z,w)=(z-w)a_{0}(z,w)$; we have
$\a(z,w)=-(z-w)\tilde a_{0}(z,w)$ and
$$
[z,\bar\rho](f)=\res_{w}\a(z,w)(a_{0}+\tilde
a_{0})(z,w)f(w)\omega_{w}.
$$
Now $a_{0}+\tilde
a_{0}=\delta(z/w)/r_{0}(z)$ so 
$$
[z,\bar\rho](f)(z)=\a(z,z)f(z).
$$
But $a_{0}(z,w)=(\sum_{i\le -N}
(z/w)^{i}+\sum_{i>-N}z^{i}\la_{i}(w))/r_{0}(z)$, with $\la_{i}\in\La$, so
that
$$
\eqalign{
(z-w)a_{0}(z,w) & =-(w/r_{0}(z))(z/w)^{-N}+(z-w) \sum_{i\ge
-N}z^{i}\la_{i}(w)/r_{0}(z) \cr &
\in
-z/r_{0}(z)
+(z-w)\CC((z,w)),
}$$
so $\a(z,z)=-z/r_{0}(z)$; so that $[\bar\rho,
z]=[\pr, z]$; so we have, $\bar\rho=\pr+$function. Since $\La$
is isotropic, $\bar\rho(1)=0$; this shows $\bar\rho=\pr$. 

Let us consider now $\pr_{z}a_{0}(z,w)-a_{0}(z,w)^{2}$; this expression
belongs to $R_{z}\otimes k_{w}$ and the endomorphism of $R$ it defines
is zero, so it belongs to $R\otimes R$. We have shown: 

\proclaim{Proposition 1} There exists $\gamma\in R\otimes R$, such that
$\pr_{z}a_{0}(z,w)=a_{0}(z,w)^{2}+\gamma(z,w)$. 
\endproclaim 

\medskip 
\noindent {\it 3.2. The kernel of $\pr^{n}$, $n\ge 0$.}

In this section we will study the expressions $\sum \pr^{k}e^{i}\otimes
e_{i}$. For $k=0$, this expression is equal to $a_{0}$; for $k=1$, it is
equal to $\pr_{z}a_{0}=a_{0}^{2}+\gamma$; for $k=2$, it is equal to 
$\pr_{z}(a_{0}^{2}+\gamma)=2a_{0}(a_{0}^{2}+\gamma)+\pr_{z}\gamma=2a_{0}^{3}
+2a_{0}\gamma+\pr_{z}\gamma$. More generally, we have: 

\proclaim{Proposition 2} Let $(P_{k}^{(n)})_{k\in\ZZ,n\ge 0}$ be the system
of polynomials in $\CC[\gamma_{0},\gamma_{1},...]$ defined by
$P^{(0)}_{<0}=0$, 
$P_{1}^{(0)}=1$, $P_{>1}^{(0)}=0$, and $P_{k}^{(n+1)}=DP_{k}^{(n)}
+(k-1)P_{k-1}^{(n)}+\gamma_{0}(k+1)P_{k+1}^{(n)}$, for $n\ge 0$, 
where $D=\sum_{i\ge
0}\gamma_{i+1} \pr/\pr\gamma_{i}$. Then 
$$
\sum\pr^{n}e^{i}\otimes e_{i}=\sum_{k\ge
0}P_{k}^{(n)}(\gamma, \pr_{z}\gamma,...)a_{0}^{k},
\leqno(3.3)
$$
$$
\sum e^{i}\otimes\pr^{n}e_{i}=(-1)^{n}\sum_{k\ge
0}P_{k}^{(n)}((-1)^{i}\pr_{w}^{i}\tilde \gamma)a_{0}^{k}.
\leqno(3.4)
$$
\endproclaim

The proof is a simple induction. Then, we have 
$\sum(\sum_{n\ge1}\hbar^{n}\pr^{n-1}/n!)e^{i}\otimes e_{i}=\sum_{k\ge
0}(\sum_{n\ge1}{\hbar^{n}\over n!}P_{k}^{(n-1)}(\gamma,\pr_{z}
\gamma,...))a_{0}^{k}$; let us set
$$
u_{k}(\hbar,\gamma_{0},\gamma_{1},...)=\sum_{n\ge
1}{\hbar^{n}\over n!}
P_{k}^{(n-1)}(\g_{0},\g_{1},...).\leqno(3.5) 
$$
We have, with $t$ an auxiliary variable, and with
$u(\hbar,t,\gamma_{i})=\sum_{k\ge 0}t^{k}u_{k}$,
the equations
$$
{{\pr u}\over{\pr\hbar}}=t+Du+(t^{2}+\gamma_{0}){{\pr u}\over{\pr t}}, 
\quad u_{|\hbar=0}=0;
\leqno(3.6) 
$$
by Cauchy's principle, they determine uniquely $u$. 
Let $\phi,\psi\in\hbar\CC[\g_{i}][[\hbar]]$ be the solutions to 
$$
{{\pr\psi}\over{\pr\hbar}}=D\psi-1-\g_{0}\psi^{2}, \quad
{{\pr\phi}\over{\pr\hbar}}=D\phi-\g_{0}\psi, \leqno(3.7)
$$
then the expansions of $\phi$ and $\psi$ are $\psi=-\hbar+...$,
$\phi=\hbar^{2}\g_{0}+...$, and $\phi$, $\psi$ have the properties 
$$
\phi(-\hbar,(-1)^{i}\g_{i})=\phi(\hbar,\g_{i}), \quad
\psi(-\hbar,(-1)^{i}\g_{i})=-\psi(\hbar,\g_{i}). \leqno(3.8)
$$
Moreover, $\phi-\ln(1+a_{0}\psi)$ satisfies (3.6); this identifies this
function with $u$. 

We conclude from this the first part of 

\proclaim{Proposition 3} With $\phi$ and $\psi$ the solutions of (3.7),
we have 
$$
\sum_{i}{{q^{\pr}-1}\over\pr}e^{i}\otimes
e_{i}=\phi(\hbar,\pr_{z}^{i}\g)
-\ln(1+a_{0}\psi(\hbar,\pr_{z}^{i}\g)),\leqno(3.9)
$$
$$
\sum_{i}{{1-q^{-\pr}}\over\pr}e^{i}\otimes
e_{i}=-\phi(-\hbar,\pr_{z}^{i}\g)
+\ln(1+a_{0}\psi(-\hbar,\pr_{z}^{i}\g)),\leqno(3.10)
$$
$$
\sum_{i}e^{i}\otimes{{q^{\pr}-1}\over\pr}
e_{i}=-\phi(\hbar,\pr_{w}^{i}\tilde\g)
+\ln(1-a_{0}\psi(\hbar,\pr_{w}^{i}\tilde\g)),\leqno(3.11)
$$
$$
\sum_{i}e^{i}\otimes{{q^{-\pr}-1}\over\pr}
e_{i}=-\phi(-\hbar,\pr_{w}^{i}\tilde\g)
+\ln(1-a_{0}\psi(-\hbar,\pr_{w}^{i}\tilde\g)).\leqno(3.12)
$$
\endproclaim
The last part is proved using the last part of prop. 2. These results
imply the following statement: 
$$
\eqalign{
\phi(\hbar,\pr_{z}^{i}\g)-\ln(1+a_{0}\psi(\hbar,\pr_{z}^{i}\g))
-\phi(-\hbar,\pr_{w}^{i}\tilde\g) & 
+\ln(1-a_{0}\psi(-\hbar,\pr_{w}^{i}\tilde
\g)) \cr & 
={{f(\pr_{z})-f(-\pr_{w})}\over{\pr_{z}+\pr_{w}}}(\g-\tilde\g), 
}
\leqno(3.13)
$$
with $f(x)={{q^{x}-1}\over x}$; this is proven by noting that the
l.h.s. belongs to $R\otimes R$, and that
$(\pr_{z}+\pr_{w})a_{0}=\g-\tilde\g$. In particular, we have
$$
\phi(\hbar,\pr_{z}^{i}\g)-\phi(-\hbar,\pr_{z}^{i}\g)
+\ln(-\psi(-\hbar,\pr_{z}^{i}\g)/\psi(\hbar,\pr_{z}^{i}\g))
={{f(\pr_{z})-f(-\pr_{w})}\over{\pr_{z}+\pr_{w}}}(\g-\tilde\g)
+\nu_{0},
\leqno(3.14)
$$
$\nu_{0}\in R\otimes R$, vanishing on the diagonal.

\proclaim{Proposition 4} Recall that
$T={{q^{\pr}-q^{-\pr}}\over{2\hbar\pr}}$; let us 
set 
$$
\psi_{0}={1\over{2\hbar}}(\phi(\hbar,\pr_{z}^{i}\g)
-\phi(-\hbar,\pr_{z}^{i}\g)), \quad
\psi_{+}(\g_{i})=\psi(-\hbar,\g_{i}), 
\quad
\psi_{-}(\g_{i})=\psi(\hbar,\g_{i}),
$$
then
$$
\sum Te^{i}\otimes e_{i}
=\psi_{0}(\gamma,\pr_{z}\gamma,...)+
{1\over{2\hbar}}\ln{{1+a_{0}\psi_{+}(\g,\pr_{z}\g,...)}\over
{1+a_{0}\psi_{-}(\g,\pr_{z}\g,...)}},
\leqno(3.15) 
$$
and 
$$
\sum e^{i}\otimes Te_{i}=-\psi_{0}(\pr_{w}^{i}\tilde\g)
+{1\over{2\hbar}}\ln{{1-a_{0}\psi_{-}(\pr_{w}^{i}\tilde\g)}
\over{1-a_{0}\psi_{+}(\pr_{w}^{i}\tilde\g)}}.
\leqno(3.16) 
$$
\endproclaim

%

From (3.13) follows
$$
\eqalign{
\psi_{0}(\pr_{z}^{i}\g)+{1\over{2\hbar}}\ln{{1+a_{0}\psi_{+}(\pr_{z}^{i}\g)}
\over{1+a_{0}\psi_{-}(\pr_{z}^{i}\g)}}
+\psi_{0}(\pr_{w}^{i}\tilde\g) & -{1\over{2\hbar}}\ln{{1-a_{0}\psi_{-}
(\pr_{w}^{i}\tilde\g)}\over{1-a_{0}\psi_{+}(\pr_{w}^{i}
\tilde\g)}} \cr
 & ={{T_{z}-T_{w}}\over{\pr_{z}+\pr_{w}}}(\g-\tilde\g). 
}\leqno(3.17)$$
Remark that 
$$
\sum Te^{i}\otimes e_{i}-e^{i}\otimes
Te_{i}={{T_{z}-T_{w}}\over{\pr_{z}+\pr_{w}}}(\g-\tilde\g); 
\leqno(3.18) 
$$
recall that $\tau=\sum Ue_{i}\otimes e^{i}$ satisfies
$$
\tau+\tilde\tau={{T_{z}-T_{w}}\over{\pr_{z}+\pr_{w}}}(\g-\tilde\g). 
\leqno(3.19) 
$$

Let us precise now the $e-e$ and $f-f$ relations. Let 
$$
A=\sum e^{i}\otimes (T+U)e_{i}, 
\leqno(3.20) 
$$
then $A=
\sum Te^{i}\otimes e_{i}+e^{i}\otimes [(Te_{i})_{R}+Ue_{i}]$; so 
$$
\eqalign{
A & =-\tau+\psi_{0}(\g,\pr_{z}\g,...)+{1\over{2\hbar}}\ln{{1+
a_{0}\psi_{+}(\g,\pr_{z}\g,...)}
\over{1+a_{0}\psi_{-}(\g,\pr_{z}\g,...)}}
\cr
& =
\tilde\tau-\psi_{0}(\tilde\g,\pr_{w}\tilde\g,...)+{1\over{2\hbar}}\ln{{1-
a_{0}\psi_{-}(\tilde\g,\pr_{w}\tilde\g,...)}
\over{1-a_{0}\psi_{+}(\tilde\g,\pr_{w}\tilde\g,...)}}.
}\leqno(3.21) $$
The relations are then written 
$$
\eqalign{
e^{2\hbar\psi_{0}(\g,\pr_{z}\g,...)}
[z-w+\a(z,w) & \psi_{+}(\g,\pr_{z}\g,...)]
e(z)e(w)= \cr & 
e^{2\hbar\tau(z,w)}[z-w+\a(z,w)\psi_{-}(\g,\pr_{z}\g,...)]e(w)e(z)}
\leqno(3.22) 
$$
$$\eqalign{
q^{K(\pr_{z}+\pr_{w})} & \{e^{2\hbar\tau(z,w)}[z-w+\a(z,w)\psi_{-}
(\g,\pr_{z}\g,...)]\}
f(z)f(w)
=\cr
&
q^{K(\pr_{z}+\pr_{w})}\{e^{2\hbar\psi_{0}(\g,\pr_{z}\g,...)}
[z-w+\a(z,w) \psi_{+}(\g,\pr_{z}\g,...)]\}
f(w)f(z)
; }
\leqno(3.23) $$
recall that $\a(z,w)=(z-w)a_{0}(z,w)$ belongs to
$(\oplus_{i}k_{z_{i}})^{\hat\otimes 2}$. 

Let us show that the $e-e$ and $f-f$ relations define a flat
deformation of the symmetric algebras in the $e(\eps)$ and $f(\eps)$,
$\eps\in\oplus_{i}k_{x_{i}}$. This statement is equivalent to 
$$
\eqalign{
 e^{2\hbar[-\psi_{0}(\g,\pr_{z}\g,...)+\tau]} & {{z-w+\a\psi_{-}(\g,
\pr_{z}\g,...)}\over{z-w+\a\psi_{+}(\g,\pr_{z}\g,...)}}\cdot \cr 
&
\Big\{ 
e^{2\hbar[-\psi_{0}(\g,\pr_{z}\g,...)+\tau]}{{z-w+\a\psi_{-}(\g,
\pr_{z}\g,...)}\over{z-w+\a\psi_{+}(\g,\pr_{z}\g,...)}}
\Big\}^{\widetilde{\ 	}}
=1, 
\cr}
$$
and this is in turn written 
$$
e^{2\hbar(-\psi_{0}-\tilde\psi_{0})}e^{2\hbar
{{T_{z}-T_{w}}\over{\pr_{z}+\pr_{w}}}(\g-\tilde\g)}
{{1+a_{0}\psi_{-}(\g,\pr_{z}\g,...)}\over{1+a_{0}\psi_{+}(\g,\pr_{z}\g,...)}}
{{1-a_{0}
\psi_{-}(\tilde\g,\pr_{w}\tilde\g,...)}\over{1+a_{0}\psi_{+}(\tilde\g,
\pr_{w}\tilde\g,...)}}=1,
$$
which amounts to the statement (3.17) above. 

To summarize, we have: 

\proclaim{Theorem 5} Let $\tau\in R\otimes R[[\hbar]]$ satisfy (2.1). 
The algebra $U_{\hbar,\La,\tau}\hat{\bf g}$ defined by generators
$K$, $D$, $h^{+}(r)$, $h^{-}(\la)$, $e(\eps)$, $f(\eps)$, 
$\la\in\La$,
$r\in R$, $\eps\in\oplus_{i}k_{x_{i}}$, subject to relations
(1.1-4) organized in generating series (1.5),
subject to relations 
$$
[h^{+}(z),e(w)]=2(\sum_{i}e_{i}\otimes e^{i})e(w), 
$$
$$
[h^{-}(z),e(w)]=2\big(\sum_{i}q^{-K\pr}e^{i}\otimes 
(T+U)e_{i}\big)e(w), 
$$
$$
[h^{+}(z),f(w)]=-2(\sum_{i}e_{i}\otimes q^{K\pr}e^{i})f(w), 
$$
$$
[h^{-}(z),f(w)]=-2\big(\sum_{i}q^{K\pr}e^{i}
\otimes 
q^{K\pr}(T+U)e_{i}\big)f(w), 
$$
$$
[h^{+}(z), h^{+}(w)]=0, \quad
[h^{+}(z), h^{-}(w)]=
{2\over\hbar} \sum_{i}e_{i}\otimes
(q^{K\pr}-q^{-K\pr})e^{i}, 
$$
$$
[h^{-}(z), h^{-}(w)]={1\over\hbar}
(q^{K\pr}\otimes q^{K\pr}-q^{-K\pr}\otimes q^{-K\pr})\sum_{i}
e^{i}\otimes (T+U)e_{i}, 
$$
$$
[e(z), f(w)]=(q^{K\pr_{w}}\delta(z/w))q^{(T+U)h^{+}(z)}
-(q^{-K\pr_{w}}\delta(z/w))q^{-h^{-}(w)},
$$
where the variables $z$ and $w$ are affected respectively to the first
and second factor; $K$ is central, (3.22), (3.23); (2.21), (2.19),
(2.20), (2.24), (2.25);  
with coproduct defined by (2.3), (2.8), (2.9), (2.22), counit defined to be
zero on all generators, and skew antipode defined by (2.12), (2.13),
(2.23), is
a Hopf algebra, quantizing the Manin triple of section 1.4. 
\endproclaim

\section{4. Dependence in $\tau$ and $\La$.}{}

\medskip
\noindent{\it 1. Dependence in $\tau$. \rm}

Let us study the dependence of the algebra $U_{\hbar,\La,\tau}\hat{\bf
g}$ defined in thm. 5, with respect to $\tau$. Let $\tau'=\tau+\u$, $\u\in
\wedge^{2}R[[\hbar]]$. 
Let us denote with a prime all quantities corresponding to
the algebra $U_{\hbar,\La, \tau'}\hat{\bf g}$. Let us denote by $u:\La\to 
R[[\hbar]]$
the linear map defined by 
$u(\la)=\langle \u,1\otimes\la \rangle_{\oplus_{i}k_{x_{i}}}$. 
We have 
$$
u=U'-U, \quad \tau'-\tau=\u=\sum_{i}(U'-U)e_{i}\otimes e^{i}. \leqno(4.1)
$$
Then: 

\proclaim{Proposition 6} The formulae 
$$
i(e'(z))=e^{{1\over 2}\hbar 
u(h^{+})(z)}e(z),
\quad
i(f'(z))=f(z)e^{{\hbar\over 2}(q^{K\pr_{z}})u(h^{+})(z)},\leqno(4.2)
$$
$$
i(h^{+\prime}(z))
=h^{+}(z), \quad 
i(h^{-\prime}(z))
=h^{-}(z)-({{q^{K\pr}+q^{-K\pr}} \over 2}uh^{+})(z),\leqno(4.3)
$$
$i(K')=
K$, $i(D')=
D$, define an algebra
isomorphism $i:U_{\hbar,\La, \tau'}\hat{\bf g}\to
U_{\hbar,\La, \tau}\hat{\bf g}$. 
Moreover, we have
$$
\Delta(i(x))=\Ad \exp({\hbar\over 4}(h^{+}\otimes h^{+})\u)\{(i\otimes
i)\Delta'(x)\}, \quad \forall x\in U_{\hbar,\La, \tau'}\hat{\bf g},
\leqno(4.4)
$$
so that both Hopf algebra structures are isomorphic up to a twisting
operation. 
\endproclaim
\medskip
\noindent 
\bf Proof. \rm
$i$ is well-defined, because $uh^{+}(z)$ is expressed as $\sum_{i\ge
0}\hbar^{i}\sum_{j}h^{+}(r^{(i)}_{j})r^{\prime(i)}_{j}(z)$, $r^{(i)}_{j}$,
$r^{\prime(i)}_{j}\in R$, the sums
$\sum_{j}h^{+}(r^{(i)}_{j})r^{\prime(i)}_{j}(z)$ being finite. 
To prove e.g. 
that $i$ is an algebra morphism, we make use (while checking the
$e-f$ relations) of the following sequence of identities:
$$
\eqalign{
& ( q^{-K\pr_{w}} \delta(z  /w)) 
e^{{\hbar\over 2}(uh^{+})(z)} 
q^{-h^{-}(w)}
e^{{\hbar\over 2} (q^{K\pr}uh^{+})(w)} 
=
\cr
&
=(q^{-K\pr_{w}}\delta(z/w)) 
e^{{\hbar\over 2}(uh^{+})(z)}
q^{-h^{-}(w)}
e^{{\hbar\over 2} (q^{2K\pr}uh^{+})(z)} 
\cr
&=
(q^{-K\pr_{w}}\delta(z/w))e^{-{\hbar\over 2}[((1-q^{2K\pr})uh^{+})(z),
h^{-}(w)]}
q^{-h^{-}(w)+{1\over 2}( ( 1+q^{-2K\pr}) uh^{+})(z)}
\cr
&=
(q^{-K\pr_{w}}\delta(z/w))e^{-{\hbar\over 2}
\{(1-q^{2K\pr})
\otimes {2\over\hbar}(q^{K\pr}-q^{-K\pr})\}
(\sum ue_{i}\otimes e^{i})}
q^{-h^{-}(w)+{1\over 2}(  ( 1+q^{-2K\pr})  uh^{+})(z)}
\cr
&=
q^{-K\pr}\{\delta(z/w)e^{(q^{2K\pr}-1)\otimes(q^{2K\pr}-1)\u}\}
q^{-h^{-}(w)+{1\over 2}((q^{K\pr}+q^{-K\pr})
uh^{+})(w)}
\cr
&=
(q^{-K\pr_{w}}\delta(z/w))q^{-h^{-}(w)+{1\over 2}((q^{K\pr}+q^{-K\pr})
uh^{+})(w)}
\cr}
$$
(the first identity follows from $\delta(z/w)f(z)=\delta(z/w)f(w)$, the
second from $e^{a}e^{b}=e^{b}e^{a}e^{[a,b]}$ if $[a,b]$ is scalar, the
last one from the fact that $\u$ is antisymmetric, so that 
$\{(q^{2K\pr}-1)\otimes(q^{2K\pr}-1)\}\u$ vanishes on the diagonal). 
The other identities are easily checked. While checking the twisting
identity for $f'(z)$, we use also the fact that 
$$
\eqalign{
[((q^{K\pr}-q^{-K\pr})uh^{+})(z),h^{-}(z)]&=
\sum_{i}(q^{K\pr}-q^{-\pr})\u_{i}(z)[h^{+}(\u'_{i}),h^{-}(z)]
\cr
&=\sum_{i}(q^{K\pr}-q^{-K\pr})
\u_{i}(z){2\over\hbar}(q^{K\pr}-q^{-K\pr})\u'_{i}(z)
\cr
&=0
,\cr
}
$$
with $\u=\sum_{i}\u_{i}\otimes\u'_{i}$, because $\{(q^{K\pr}-q^{-K\pr})
\otimes(q^{K\pr}-q^{-K\pr})\}\u\in\wedge^{2}R[[\hbar]]$. 
\bull

\proclaim{Proposition 7} The formulae 
$$
i'(e'(z))=e(z)e^{{1\over 2}\hbar 
u(h^{+})(z)},
\quad
i'(f'(z))=
e^{{\hbar\over 2}(q^{K\pr_{z}})u(h^{+})(z)}f(z),\leqno(4.5)
$$
$$
i'(h^{+\prime}(z))
=h^{+}(z), \quad 
i'(h^{-\prime}(z))
=h^{-}(z)-({{q^{K\pr}+q^{-K\pr}} \over 2}uh^{+})(z),\leqno(4.6)
$$
$i'(K')=
K$, $i'(D')=
D$, also define an algebra
isomorphism $i':U_{\hbar,\La, \tau'}\hat{\bf g}\to
U_{\hbar,\La, \tau}\hat{\bf g}$, satisfying  
$$
\Delta(i'(x))=\Ad \exp({\hbar\over 4}(h^{+}\otimes h^{+})\u)\{(i'\otimes
i')\Delta'(x)\}, \quad \forall x\in U_{\hbar,\tau'}\hat{\bf g}.
\leqno(4.7)
$$
It follows, that $i^{\prime -1}\circ i$ is a Hopf algebra automorphism
of $U_{\hbar,\tau'}\hat{\bf g}$.
\endproclaim

\medskip
\it 2. Dependence in $\La$. \rm 

Let $\La$ and $\bar\La$ be two Lagrangian supplementaries to $R$. Then we
have, $\bar\La=(1+r)\La$, with $r:\La\to R$, given by 
$$
r(\la)=\langle r_{0}, 1\otimes\la\rangle, \quad r_{0}\in\wedge^{2}R.
\leqno(4.8)
$$
Dual bases for $R$ and $\bar\La$ are then $(e^{i})$ and $(\bar e_{i})$, with
$\bar e_{i}=(1+r)e_{i}$. 
Let us set 
$$
\bar\tau=\sum_{i} U\bar e_{i}\otimes e^{i}=\tau-\sum_{i} Tre_{i}\otimes e^{i};
\leqno(4.9)
$$
we have then 
$$
\sum_{i} (T+\bar U)\bar e_{i}\otimes e^{i}=\sum_{i} (T+U)e_{i}\otimes
e^{i}.
\leqno(4.10)
$$
Let us consider the Hopf algebras $U_{\hbar, \La,\tau}\hat{\bf g}$,
$U_{\hbar, \bar\La,\bar\tau}\hat{\bf g}$ and let us denote with a bar the
quantities occurring in the second. 

\proclaim{Proposition 8} The mapping 
$$
j:U_{\hbar,\La,\tau}\hat{\bf g}\to U_{\hbar,\bar\La,\bar\tau}\hat{\bf g}
$$
defined by $j(\bar e(z))=e(z)$, $j(\bar f(z))=f(z)$, 
$j(\bar h^{+}(e^{i}))=h^{+}(e^{i})$, 
$j(\bar h^{-}(\bar e_{i}))=h^{-}(e_{i})$, $j(\bar D)=D$, $j(\bar K)=K$,
defines a Hopf algebras isomorphism between 
$U_{\hbar,\La,\tau}\hat{\bf g}$ and $U_{\hbar,\bar\La,\bar\tau}\hat{\bf g}$. 
\endproclaim

\section{5. Finite dimensional representations.}{}

Let us fix $\La$ and $\tau$, and denote by $U_{\hbar}\hat{\bf
g}_{|K=0,{\rm no}D}$ the algebra defined in thm. 5, without
generator $D$ and with $K$ specialized to zero. We construct a morphism
of algebras 
$$
\pi: U_{\hbar}\hat{\bf g}_{|K=0,{\rm no}D}\to
\End(\CC^{2})\otimes (\oplus_{i}k_{x_{i}})[[\hbar]],
\leqno(5.1)
$$
as follows: let us denote 
by $\zeta=(\zeta_{i})$ the system of coordinates $(z_{i})$, occurring
in the r.h.s.; we define 
$$
\pi(h^{+}(r))=r(\zeta)h+\rho^{+}(r)(\zeta){\rm Id}_{\CC^{2}},
\pi(h^{-}(\la))=(T+U)(\la)(\zeta)h+\rho^{-}(\la){\rm Id}_{\CC^{2}}, 
\leqno(5.2) 
$$
$$
\pi(e(z))=F(\zeta)\delta(z/\zeta)e, \quad \pi(f(z))=\delta(z/\zeta)f, 
\leqno(5.3) 
$$
where the $e,f,h$ occurring in the r.h.s. are the matrices with nonzero
coefficients $e_{12}=f_{21}=h_{11}=-h_{22}=1$, 
and 
$$
\rho^{+}:R\to R[[\hbar]], \rho^{-}: \La\to
(\oplus_{i}k_{x_{i}})[[\hbar]],
F(\zeta)\in (\oplus_{i}k_{x_{i}})[[\hbar]]
$$ 
are subject 
to the following conditions: recall that 
$A(\zeta,z)=\sum e^{i}(\zeta)(T+U)(e_{i})(z)$, and let 
$$
\b(\zeta,z)=\sum
\rho^{+}(e^{i})(\zeta)(T+U)(e_{i})(z),
\g(\zeta,z)=\sum \rho^{-}(e_{i})(\zeta)e^{i}(z); 
\leqno(5.4) 
$$
then 
$$
q^{A+\b}-q^{-\tilde A-\g}=F(z)\delta(z/\zeta), 
q^{-A+\b}-q^{\tilde A-\g}=-F(z)\delta(z/\zeta).
\leqno(5.5) 
$$
Let $F(z)\delta(z/\zeta)=e^{\sigma}\hbar(a_{0}+\tilde a_{0})(z,\zeta)$,
with $\sigma\in \hbar R\otimes R[[\hbar]]$, then we have
for some $\rho_{1,2}\in \hbar^{-1}+R\otimes R[[\hbar]]$, 
$$
A={1\over {2\hbar}}\ln{{\rho_{1}+a_{0}}\over{\rho_{2}-a_{0}}}, 
\b={1\over{2\hbar}}\ln\hbar^{2}(\rho_{1}+a_{0})(\rho_{2}-a_{0})+\sigma,
\leqno(5.6)
$$
$$
\tilde A={1\over {2\hbar}}\ln{{\rho_{2}+\tilde
a_{0}}\over{\rho_{1}-\tilde a_{0}}}, 
\g=-{1\over{2\hbar}}\ln\hbar^{2}(\rho_{1}-\tilde a_{0})(\rho_{2}
+\tilde a_{0})-\sigma.\leqno(5.7)
$$

Let us determine the possible $\rho_{1,2}$ satisfying (5.6.a), (5.7.a).
Comparing (5.6.a) and the second line of (3.21), it is enough to have
$$
\ln\hbar(\rho_{1}+a_{0})=\ln(1-a_{0}\psi_{-}
(\pr^{i}_{w}\tilde\gamma))+2\hbar\la
\leqno(5.8) 
$$
$$
\ln\hbar(\rho_{2}-a_{0})=\ln(1-a_{0}\psi_{+}(\pr^{i}_{w}\tilde\gamma))
+2\hbar\bar\la
\leqno(5.9) 
$$
with $\la, \bar\la\in R\otimes R[[\hbar]]$,
$\la-\bar\la=\tilde\tau-\psi_{0}(\partial_{w}^{i}\tilde\gamma)$; and
comparing (5.7.a) and the first line of (3.21), it is enough to have
$$
\ln\hbar(\tilde\rho_{2}+a_{0})=\ln(1+a_{0}\psi_{+}(\pr_{z}^{i}\gamma))+2\hbar
\mu
\leqno(5.10) 
$$
$$
\ln\hbar(\tilde\rho_{1}-a_{0})=\ln(1+a_{0}\psi_{-}(\pr_{z}^{i}\gamma))
+2\hbar\bar\mu
\leqno(5.11) 
$$
with $\mu,\bar\mu\in R\otimes R$, and
$\mu-\bar\mu=\psi_{0}(\pr_{z}^{i}\gamma)-\tau$. 


(5.9-11) are
equivalent to the fact that for certain $\nu,\nu'\in R\otimes R$, equal
to $0$ on the diagonal,  
$$
\bar\la=\tilde\mu=\nu+{1\over{2\hbar}}
\ln(\hbar/\psi_{+}(\pr_{w}^{i}\tilde\g)), \rho_{2}
={e^{\nu}\over{\psi_{+}(\pr_{w}^{i}\tilde\g)}}+a_{0}(1-e^{\nu}), 
\leqno(5.12)
$$
and
$$
\la=\tilde{\bar\mu}=\nu'+{1\over{2\hbar}}
\ln(-\hbar/\psi_{-}(\pr_{w}^{i}\tilde\g)),
\rho_{1}=-{e^{\nu'}\over{\psi_{-}(\pr_{w}^{i}\tilde\g)}}+a_{0}(e^{\nu'}-1),
\leqno(5.13)
$$
with the conditions on $\nu$ and $\nu'$
$$
\nu'-\nu+{1\over{2\hbar}}\ln(\psi_{+}(\pr_{w}^{i}\tilde\gamma)
/-\psi_{-}(\pr_{w}^{i}\tilde\gamma))
=\tilde\tau-\psi_{0}(\pr_{w}^{i}\tilde\gamma).
\leqno(5.14)
$$

Let us see now, how $\rho^{\pm}$ can be deduced from these equalities.
The conditions on them are 
$$
\ln(1-a_{0}\psi_{-}
(\pr^{i}_{w}\tilde\gamma))+2\hbar\la
=\hbar(1+\rho^{+})(e^{i})\otimes (T+U)(e_{i})-\hbar\sigma,
\leqno(5.15)$$
$$
\ln(1-a_{0}\psi_{+}(\pr^{i}_{w}\tilde\gamma))
+2\hbar\bar\la
=\hbar(\rho^{+}-1)(e^{i})\otimes (T+U)(e_{i})-\hbar\sigma,
\leqno(5.16)
$$
$$
\ln(1+a_{0}\psi_{+}(\pr_{z}^{i}\gamma))+2\hbar
\mu
=\hbar\{e^{i}\otimes (T+U)(e_{i})-e^{i}\otimes \rho^{-}(e_{i})\}
-\hbar\tilde\sigma,
\leqno(5.17)
$$
$$
\ln(1+a_{0}\psi_{-}(\pr_{z}^{i}\gamma))
+2\hbar\bar\mu
=-\hbar\{e^{i}\otimes (T+U)(e_{i})+e^{i}\otimes \rho^{-}(e_{i})\}
-\hbar\tilde\sigma. 
\leqno(5.18)
$$

Let $T_{+}$, $T_{-}$ be the endomorphisms of $R$, defined by 
$$
T_{\pm}(r)=\langle \ln(1-a_{0}\psi_{\mp}(\pr_{w}^{i}\tilde\gamma)),
1\otimes r\rangle; 
\leqno(5.19)
$$
we have $T_{+}={{1-q^{-\pr}}\over{\pr}}$,
$T_{-}={{1-q^{\pr}}\over{\pr}}$. Since
$T_{\pm}=\hbar(\rho^{+}\pm 1)T$, (recall that
$T={{q^{\pr}-q^{-\pr}}\over
{2\hbar\pr}}$),  
$$
\rho^{+}={{1-q^{\pr}}\over{1+q^{\pr}}}. 
\leqno(5.20)
$$
Due to (5.12) (resp. (5.13)), (5.15) and (5.16) (resp. (5.17) and
(5.18)) are
equivalent. (5.14) can be solved by posing
$$
\tau={1\over{2\hbar}}{{f(\pr_{z})-f(-\pr_{w})}\over{\pr_{z}+\pr_{w}}}
(\g-\tilde\g), \quad \nu'=0,\nu=\tilde\nu_{0}. \leqno(5.21)
$$
(5.15) then gives us 
$$
\eqalign{
\hbar\sigma={1\over{\pr_{z}+\pr_{w}}}[{1\over{1+q^{\pr_{z}}}}(f(-\pr_{z})
+f(-\pr_{w})) & -f(\pr_{w})](\g-\tilde\g) \cr & 
+\phi(\hbar,\pr_{w}^{i}\tilde\g)
-\ln(-\hbar/\psi_{-}(\pr_{w}^{i}\tilde\g)),
}
\leqno(5.22)
$$
and (5.18) gives then
$$
\eqalign{
\sum & e^{i}\otimes\rho^{-}(e_{i})
=\sum e^{i}\otimes
{{(q^{\pr}-1)(1-q^{-\pr})}\over{2\hbar\pr}}e_{i}-{1\over\hbar}
[\phi(-\hbar,\pr_{z}^{i}\g)+\phi(\hbar,\pr_{z}^{i}\g)]
\cr
&
-{1\over\hbar}{1\over{\pr_{z}+\pr_{w}}}\big[{1\over
2}(f(\pr_{w})-f(-\pr_{z}))-{1\over{1+q^{\pr_{w}}}}(f(-\pr_{w})+f(-\pr_{z}))
+f(\pr_{z})\big](\g-\tilde\g)
}
$$
and so 
$$
\eqalign{
\rho^{-}(\la)= & {{(q^{\pr}-1)(1-q^{-\pr})}\over{2\hbar\pr}}\la
-{1\over\hbar}\langle
[\phi(-\hbar,\pr_{z}^{i}\g)+\phi(\hbar,\pr_{z}^{i}\g)]
+{1\over{\pr_{z}+\pr_{w}}}\big[{1\over
2}\big(f(\pr_{w}) \cr & -f(-\pr_{z})\big) 
-{1\over{1+q^{\pr_{w}}}}(f(-\pr_{w})+f(-\pr_{z}))
+f(\pr_{z}) \big](\g-\tilde\g), \la\otimes 1 \rangle
}
\leqno(5.23)
$$
for $\la\in\La$. So we have:

\proclaim{Proposition 9} The formulae (5.2), (5.3) define a morphism of
algebras 
$$
\pi: U_{\hbar}\hat{\bf g}_{|K=0,{\rm no}D}\to
\End(\CC^{2})\otimes (\oplus_{i}k_{x_{i}})[[\hbar]],
$$
provided $\tau$ is chosen according to (5.21),
with $\rho^{\pm}$, $\sigma$ given by (5.20), (5.23) and (5.22). 
\endproclaim

Let us indicate how the formulae giving $\rho^{\pm}$ and $\sigma$ would
be altered in the case of an arbitrary $\tau$ (satisfying (2.1)). Let us
denote with an exponent $ ^{(0)}$ the quantities implied in prop. 9. 
The general form of a solution of (2.1) is $\tau=\tau^{(0)}+\a$,
$\a\in\wedge^{2}R[[\hbar]]$; we have then
$\sigma=\sigma^{(0)}-((1+\rho^{+})\otimes 1)\a$, $\rho^{+}=\rho^{+(0)}$,
$\rho^{-}(\la)=
\rho^{-(0)}(\la)-\rho^{+}(\langle \la\otimes1,\a\rangle)$.

\section{}{6. Examples.}

\medskip
\noindent{\it 1. Trigonometric case.} 

Let $X=\CC P^{1}$, let $z$ be a coordinate on $X$, and let
$\omega=dz/z$. The set of marked points is $\{0,\infty\}$. Let us pose
$$
\La=\{(\la_{0},\la_{\infty})\in\CC[[z]]\times\CC[[z^{-1}]] |
\la_{0}(0)+\la_{\infty}(\infty)=0 \}. 
$$
Dual bases for $R$ and $\La$ are $e^{i}=z^{i}$ for $i\in\ZZ$, and
$e_{i}=(z^{-i},0)$ for $i<0$, $-(0,z^{-i})$ for $i>0$, ${1\over
2}(1,-1)$ for $i=0$. We compute then 
$$
\sum_{i}Te^{i}\otimes e_{i}={1\over{2\hbar}}(\ln{{qz-w}\over{z-qw}},0)
-{1\over{2\hbar}}(0,\ln{{qw-z}\over{w-qz}}), 
$$
so that we can take $U=0$; $\exp(2\hbar\sum_{i}Te^{i}\otimes
e_{i})=({{qz-w}\over{z-qw}},{{qz-w}\over{z-qw}})$ 
and the $e-e$ relation is 
$$
(z-qw)e(z)e(w)=(qz-w)e(w)e(z),
$$
as it appeared first in [4].

\medskip
\noindent{\it 2. Elliptic case.} 
 
Let $X$ be the elliptic curve $\CC/\ZZ+\tau_{0}\ZZ$; let $z$ be the
coordinate on $\CC$, and let $\omega=dz$. Let us consider the case
$x=0$. We choose $\La$ to be spanned by $z^{-1}, z, z^{2},z^{3},...$. 

We define $t=e^{2i\pi z}$, $q_{0}=e^{2i\pi\tau_{0}}$ (we assume $|q_{0}|<1$);
$$\theta(t)=\prod_{i\ge
0}(1-q_{0}^{n}t)\prod_{i<0}(1-q_{0}^{n}t^{-1}),\quad
\zeta={d\over{dz}}(\ln\theta).
$$

Let us compute the kernel of $T$. We have $a_{0}(z,w)
=\zeta(z-w)-\zeta(z)+\zeta(w)$, so that for $r\in R$,
$$
r(z)=\res_{w=0}(\zeta(z-w)-\zeta(z)+\zeta(w))r(w)dw,
$$
and 
$$
[(q^{\pr}-q^{-\pr})r](z)=\res_{w=0}(\zeta(z-w+\hbar)-\zeta(z+\hbar)-
\zeta(z-w-\hbar)+\zeta(z-\hbar))r(w)dw, 
$$
so
$$
[\pr^{-1}(q^{\pr}-q^{-\pr})r](z)=\res_{w=0}\ln{{\theta(z-w+\hbar)}
\over{\theta(z-w-\hbar)}}
{{\theta(z-\hbar)}\over{\theta(z+\hbar)}}
{{\theta(-w-\hbar)}\over{\theta(-w+\hbar)}}r(w)dw.
$$
We have then, 
$$
\sum \pr^{-1}(q^{\pr}-q^{-\pr})e^{i}\otimes e_{i}\in
\ln{{\theta(z-w+\hbar)}
\over{\theta(z-w-\hbar)}}
{{\theta(z-\hbar)}\over{\theta(z+\hbar)}}
{{\theta(-w-\hbar)}\over{\theta(-w+\hbar)}}
+R\otimes R[[\hbar]]
$$
and so
$$
\sum 
e^{i}\otimes [-\pr^{-1}(q^{\pr}-q^{-\pr})+U]e_{i}\in
\hbar+\ln{{\theta(z-w+\hbar)}
\over{\theta(z-w-\hbar)}}
{{\theta(z-\hbar)}\over{\theta(z+\hbar)}}
{{\theta(-w-\hbar)}\over{\theta(-w+\hbar)}}
+\wedge^{2}R\otimes R[[\hbar]];
$$
so that in the present case, the $e-e$ relation takes the form
$$
e(z)e(w)=e^{\hbar}{{\theta(z-w+\hbar)}\over{\theta(z-w-\hbar)}}
{{\theta(z-\hbar)}
\over{\theta(z+\hbar)}}{{\theta(-w-\hbar)}\over{\theta(-w+\hbar)}}e(w)e(z);
\leqno(6.1)
$$
this relation is analogous to the relation (7.3) occurring in [8]. 

\medskip
\noindent{\it 3. Double extensions and infinite twists of the Reyman-Semenov
triples.}

Let as above $X$ be an elliptic curve, and $X_{n}$ be the set of its
$n$-division points. We fix an isomorphism of $X_{n}$ with
$(\ZZ/n\ZZ)^{2}$, and denote by $a\mapsto I_{a}$ the projective
representation of $(\ZZ/n\ZZ)^{2}$ on $\oplus_{i\in\ZZ/n\ZZ}\CC
\epsilon_{i}$, 
defined by $I_{(1,0)}\epsilon_{i}=\zeta^{i}\epsilon_{i}$
$I_{(0,1)}\epsilon_{i}=\epsilon_{i+1}$, $\zeta$ being a primitive $n$-th
root of $1$. 

The following Manin triple was introduced in [12]. Let $k_{0}$,
$\cO_{0}$ be the local field and ring at $0\in X$, and let us define in  
${\bf g}=sl_{n}(k_{0})$ the scalar product $\langle x,y\rangle_{{\bf g}}
=\res_{0}\tr
(xy(z))dz$. Let ${\bf g}_{+}=sl_{n}(\cO_{0})$ and ${\bf g}_{-}$ be the
set of the 
expansions at $0$, of the regular maps $\sigma:X-X_{n}\to sl_{n}(\CC)$, 
such that $\sigma(x+a)=\Ad I_{a}\sigma(x)$, for $a\in X_{n}$. Then
$({\bf g},{\bf g}_{+},{\bf g}_{-})$ forms a Manin triple. 
Its quantization was treated in [15] in the $sl_{2}$ case, 
and is connected with Sklyanin algebras ([14]). 

We propose the following double extension for this triple. 
Let $\hat {\bf g}={\bf g}\oplus \CC K\oplus \CC D$, and let us denote
with an index 
$0$ the Lie bracket in ${\bf g}$. We endow $\hat {\bf g}$ with the bracket
$[x,y]=[x,y]_{0}+\res_{0}\tr(xdy)K$, $[D,x]={{dx}\over{dz}}$, for
$x,y\in {\bf g}$, $K$ is central. Let us also define a scalar product
$\langle, \rangle_{\hat {\bf g}}$ on $\hat {\bf g}$ by
$\langle x,y \rangle_{\hat {\bf g}}=\langle x,y \rangle_{{\bf g}}$ for
$x,y\in {\bf g}$, 
$\langle K,x\rangle_{\hat {\bf g}}=\langle D,x\rangle_{\hat {\bf g}}=0$
for $x\in 
{\bf g}$, $\langle K,D\rangle_{\hat {\bf g}}=2$. 
Then $\tilde {\bf g}_{+}={\bf g}_{+}\oplus \CC D$ and $\hat {\bf
g}_{-}={\bf g}_{-}\oplus\CC K$ 
are Lagrangian subalgebras of $\hat {\bf g}$, so that $(\hat {\bf
g},\tilde {\bf g}_{+},\hat{\bf g}_{-})$ forms a Manin triple. 

We also propose the following twist for this triple. 
Let ${\bf h}$ and 
${\bf n}_{\pm}$ be the diagonal and upper (resp. lower) triangular
subalgebras of $sl_{n}$, and let $\bar {\bf g}_{+}={\bf h}(\cO_{0})\oplus
{\bf n}_{+}(k_{0})\oplus\CC D$, 
$$
\bar {\bf g}_{-}=\{\sigma:X-X_{n}\to {\bf h}(\CC)|
\sigma(x+a)=\Ad I_{a}\sigma(x){\rm,\ for\ }\  a\in X_{n}
\}\oplus {\bf n}_{-}(k_{0})\oplus\CC K.
$$ 
Then $(\hat {\bf g},\bar {\bf g}_{+},\bar {\bf g}_{-})$ is a twist of
the previous Manin triple, which can be quantized according to the
techniques developed above in the case $n=2$.

\vskip 1truecm
\noindent
{\bf References}
\bigskip
\item{[1]} J. Beck, {\sl Braid group action and quantum affine
algebras,} Commun. Math. Phys. 165 (1994), 555-68. 
\item{[2]} V. Chari, A. Pressley, {\sl Quantum affine algebras,} 
Commun. Math. Phys. 142 (1991), 261-83. 
\item{[3]} J. Ding, I.B. Frenkel, {\sl Isomorphism of two
realizations of quantum affine algebras $U_{q}(\widehat{gl}_{n})$,} 
Commun. Math. Phys. 156 (1993), 277-300. 
\item{[4]} V.G. Drinfeld, {\sl A new realization of Yangians and
quantized affine algebras,} Sov. Math. Dokl. 36 (1988).  
\item{[5]} V.G. Drinfeld, {\sl Quasi-Hopf algebras,} Leningrad
Math. J. 1:6 (1990), 1419-57. 
\item{[6]} B. Enriquez, G. Felder, in preparation.
\item{[7]} P. Etingof, D. Kazhdan, {\sl Quantization of Lie bialgebras,
I,} q-alg/9506005.
\item{[8]} B.L. Feigin, E.V. Frenkel, {\sl Quantum $\cW$-algebras and
elliptic algebras,} RIMS preprint (1995), q-alg/9508009.  
\item{[9]} I.B. Frenkel, N. Jing, {\sl Vertex representations of
quantum affine algebras,} Proc. Natl. Acad. Sci. USA, 85 (1988), 9373-7. 
\item{[10]} S.M. Khoroshkin, V.N. Tolstoy, {\sl On Drinfeld's
realization of quantum affine algebras,} J. Geom. Phys. 11 (1993), 445-52. 
\item{[11]} N.Yu. Reshetikhin, M.A. Semenov-Tian-Shansky, {\sl Central
extensions of quantum current groups,} Lett. Math. Phys. 19 (1990), 133-42. 
\item{[12]} A.G. Reyman, M.A. Semenov-Tian-Shansky, {\sl Integrable
systems II, ch. 11,} Encycl. Sov. Math., 16,
``Dynamical systems, 7'', Springer-Verlag (1993), 188-225.  
\item{[13]} M.A. Semenov-Tian-Shansky, {\sl Poisson-Lie groups, quantum
duality principle, and the quantum double,} Theor. Math. Phys. 93
(1992), 1292-307.
\item{[14]} E.K. Sklyanin, {\sl Some algebraic structures connected
with the Yang-Baxter equation,} Funct. An. Appl. 16 (1982), 263-70. 
\item{[15]} D.B. Uglov, {\sl The quantum bialgebra associated with
the eight-vertex $R$-matrix,} Lett. Math. Phys. 28 (1993), 139-42. 
\medskip
\medskip\medskip
\section{}{}

B.E., V.R.: Centre de Math\'{e}matiques, URA 169 
du CNRS, Ecole Polytechnique, 91128 Palaiseau, France

V.R.: ITEP, 25, Bol. Cheremushkinskaya, 117259 Moscou, Russie. 

\bye